\documentclass[11pt]{article}
\usepackage{array,amssymb,amsmath,amsthm}
\usepackage[sort&compress,round]{natbib}
\usepackage{pstricks,pst-node,mathrsfs,graphicx}

\usepackage[left=3.25cm,right=3.25cm,top=3cm,bottom=3cm]{geometry}

\DeclareMathOperator{\Prob}{Prob}
\DeclareMathOperator{\E}{E}
\DeclareMathOperator{\trace}{trace}
\newcommand{\ND}{\mathcal{N}}
\newcommand{\RRR}{\mathbb{R}} 
\newcommand{\ind}{ \perp\hspace{-0.21cm}\perp }

\newcommand{\bbk}{\mathbb{K}}
\newcommand{\bbr}{\mathbb{R}}
\newcommand{\bfx}{\mathbf{x}}
\newcommand{\bfa}{\mathbf{a}}
\newcommand{\bbn}{\mathbb{N}}

\newcommand{\bft}{\mathbf{t}}

\newcommand{\bbc}{\mathbb{C}}

\newtheorem{definition}{Definition}
\newtheorem{theorem}[definition]{Theorem} 
\newtheorem{corollary}[definition]{Corollary}
\newtheorem{remark}[definition]{Remark}

\newtheorem{example}[definition]{Example}
\newtheorem{proposition}[definition]{Proposition}

\numberwithin{equation}{section}

\begin{document}
\renewcommand{\baselinestretch}{1.2}
\markright{
}
\markboth{\hfill{\footnotesize\rm MATHIAS DRTON AND SETH SULLIVANT}\hfill}
{\hfill {\footnotesize\rm ALGEBRAIC STATISTICAL MODELS} \hfill}
\renewcommand{\thefootnote}{}
$\ $\par
\fontsize{10.95}{14pt plus.8pt minus .6pt}\selectfont
\vspace{0.8pc}
\centerline{\large\bf ALGEBRAIC STATISTICAL MODELS}
\vspace{.4cm}
\centerline{Mathias Drton and Seth Sullivant}
\vspace{.4cm}
\centerline{\it University of Chicago and Harvard University}
\vspace{.55cm}
\fontsize{9}{11.5pt plus.8pt minus .6pt}\selectfont

\begin{quotation}
\noindent {\it Abstract:}
Many statistical models are algebraic in that they are defined in terms of
polynomial constraints, or in terms of polynomial or rational
parametrizations.  The parameter spaces of such models are typically
semi-algebraic subsets of the parameter space of a reference model with
nice properties, such as for example a regular exponential family.  This
observation leads to the definition of an `algebraic exponential
family'.  This new definition provides a unified framework for the
study of statistical models with algebraic structure.
In this paper we review the ingredients to this definition and illustrate
in examples how computational algebraic geometry can be used to solve
problems arising in statistical inference in algebraic models.
\par

\vspace{9pt}
\noindent {\it Key words and phrases:}
Algebraic statistics, computational algebraic geometry,
exponential family, maximum likelihood estimation, model invariants,
singularities.\\ 
\par
\end{quotation}\par

\fontsize{10.95}{14pt plus.8pt minus .6pt}\selectfont
\setcounter{section}{1}
\setcounter{equation}{0} 
\noindent {\bf 1. Introduction}

Algebra has seen many applications in statistics
\citep[e.g.][]{viana:2001,Diaconis}, but it is only rather
recently that computational algebraic geometry and related techniques in
commutative algebra and combinatorics have been used to study statistical
models and inference problems.  This use of computational algebraic
geometry was initiated in work on exact tests of conditional independence
hypotheses in contingency tables \citep{diaconis:1998}.  Another line of
work in experimental design led to the monograph by \cite{pistone:2001}.
`Algebraic statistics', the buzz word in the titles of this monograph and
the more recent book by \cite{ascb}, has now become the umbrella term for
statistical research involving algebraic geometry.  There has also begun to
be a sense of community among researchers working in algebraic statistics
as reflected by workshops, conferences, and summer schools.  One such
workshop, the 2005 Workshop on Algebraic Statistics and Computational
Biology held at the Clay Mathematics Institute led to the \emph{Statistica
  Sinica} theme topic, of which this article forms a part.  Other recent
work in algebraic statistics has considered contingency table analysis
\citep{dobra:2004,aoki:2005,takemura:2005}, phylogenetic tree models
\citep{eriksson:2005,sturmfels:2005,allman:2003}, maximum likelihood
estimation under multinomial sampling \citep{hosten:2005,catanese:2006},
reliability theory \citep{giglio:2004}, and Bayesian networks
\citep{garcia:2005}.  A special issue of the \emph{Journal of Symbolic
  Computation} emphasizing the algebraic side emerged following the 2003
Workshop on Computational Algebraic Statistics at the American Institute of
Mathematics.

The algebraic problems studied in algebraic statistics are of a rather
diverse nature.  At the very core of the field, however, lies the notion of
an algebraic statistical model.  While this notion has the potential of
serving as a unifying theme for algebraic statistics, there does not seem,
at present, to exist a unified definition of an algebraic statistical
model.  This lack of unity is apparent even when reading articles by the
same authors, where two papers might use two different, non-equivalent
definitions of an algebraic statistical model, for different theoretical
reasons.  The usual set-up for discussing algebraic statistical models has
involved first restricting to discrete random variables and then
considering models that are either conditional independence models or
defined parametrically with a polynomial or rational parametrization.
However, many statistical models for continuous random variables also have
an algebraic flavor, though currently there has been no posited description
of a general class of algebraic statistical models that would include
models for continuous random variables.

The main goal of this paper is to give a unifying definition of algebraic
statistical models, as well as illustrate the usefulness of the definition
in examples.
Our approach is based on the following philosophy.  Let
$\mathcal{P}=(P_\theta\mid \theta\in \Theta)$ be a statistical model with
parameter space $\Theta\subseteq\RRR^k$.  In this paper, a model such as
$\mathcal{P}$ is defined to be a family of probability distributions on
some given sample space.  (For a discussion of the notion of a statistical
model see \cite{mccullagh:2002} who proposes to refine the traditional
definition to one that ensures that the model extends in a meaningful way
under natural extensions of the sample space.)  Suppose that in model
$\mathcal{P}$ a statistical inference procedure of interest is
well-behaved.  If this is the case, then the properties of the inference
procedure in a submodel $\mathcal{P}_M=(P_\theta\mid \theta\in M)$ are
often determined by the geometry of the set $M\subseteq\Theta$.  Hence, if
the set $M$ exhibits algebraic structure, then the inference procedure can
be studied using tools from algebraic geometry.  This philosophy suggests
the following definition.  The semi-algebraic sets appearing in the
definition will be defined in Section 3.

\begin{definition}
  \label{def:algsubmodel}
  Let $\mathcal{P}=(P_\theta\mid \theta\in\Theta)$ be a ``well-behaved''
  statistical model whose parameter space $\Theta\subseteq\RRR^k$ has
  non-empty interior.  A submodel $\mathcal{P}_M=(P_\theta\mid \theta\in
  M)$ is an {\em algebraic statistical model} if there exists a
  semi-algebraic set $A \subseteq\RRR^k$ such that $M= A \cap \Theta$.
\end{definition}

Definition \ref{def:algsubmodel} is intentionally vague and the precise meaning of
the adjective ``well-behaved'' depends on the context.  For example, if
asymptotic properties of maximum likelihood estimators are of interest then
the word ``well-behaved'' could refer to models satisfying regularity
conditions guaranteeing that maximum likelihood estimators are
asymptotically normally distributed.  However, one class of statistical
models, namely regular exponential families, can be considered to be
well-behaved with respect to nearly any statistical feature of interest.

\begin{definition}
  \label{def:alg-exp-fam1}
  Let $(P_\eta\mid\eta\in N)$ be a regular exponential family of order $k$.
  The subfamily induced by the set $M\subseteq N$ is an \emph{algebraic
    exponential family} if there exists an open set $\bar N\subseteq
  \RRR^k$, a diffeomorphism $g:N\to \bar N$, and a semi-algebraic set $A
  \subseteq \RRR^k$ such that $M=g^{-1}(A \cap \bar N)$.
\end{definition}

Definition \ref{def:alg-exp-fam1} allows one to consider algebraic
structure arising after the regular exponential family is reparametrized
using the diffeomorphism $g$ (see Section 2.2 for a definition of
diffeomorphisms).  Frequently, we will make use of the mean
parametrization.  Algebraic exponential families appear to include all the
existing competing definitions of algebraic statistical models as special
cases.  Among the examples covered by Definition \ref{def:alg-exp-fam1} are
the parametric models for discrete random variables studied by \cite{ascb}
in the context of computational biology.
Other models included in the framework are conditional independence models
with or without hidden variables for discrete or jointly Gaussian random
variables.  Note that some work in algebraic statistics has focused on
discrete distributions corresponding to the boundary of the probability
simplex \citep{geiger:2006}.  These distributions can be included in an
extension of the regular exponential family corresponding to the interior
of the probability simplex; see \citet[pp.\ 154ff]{barndorff-nielsen:1978},
\citet[pp.\ 191ff]{brown:1986}, and \cite{csiszar:2005}.  Models given by
semi-algebraic subsets of the (closed) probability simplex can thus be
termed `extended algebraic exponential families'.

In the remainder of the paper we will explain and exemplify our definition
of algebraic exponential families.  We begin in Section 2 by reviewing
regular exponential families and in Example \ref{ex:stratified} we stress
the fact that submodels of regular exponential families are only
well-behaved if the local geometry of their parameter spaces is
sufficiently regular.  
In Section 3, we review some basic terminology and results on
semi-algebraic sets, which do have nice local geometric properties, and
introduce our algebraic exponential families.  We also show that other
natural formulations of an algebraic statistical model in the discrete case
fall under this description and illustrate the generality using jointly
normal random variables.  We then illustrate how problems arising in
statistical inference in algebraic models can be addressed using
computational algebraic geometry.  Concretely, we discuss in Section 4 how
so-called model invariants reveal aspects of the geometry of an algebraic
statistical model that are connected to properties of statistical inference
procedures such as likelihood ratio tests.  As a second problem of a
somewhat different flavour we show in Section 5 how systems of polynomial
equations arising from likelihood equations can be solved algebraically.\\
\par

\setcounter{section}{2}
\setcounter{equation}{0} 
\noindent {\bf 2. Regular exponential families}

Consider a sample space $\mathcal{X}$ with $\sigma$-algebra $\mathcal{A}$
on which is defined a $\sigma$-finite measure $\nu$.  Let $T:\mathcal{X}\to
\RRR^k$ be a statistic, i.e., a measurable map.  Define the
{\em natural parameter space\/}
\[
N = \left\{\eta\in\RRR^k \;:\; \int_{\mathcal{X}} e^{\eta^t T(x)} d\nu(x)
  <\infty\right\}.
\]
For $\eta\in N$, we can define a probability density $p_\eta$ on
$\mathcal{X}$ as
\[
p_\eta(x) = e^{\eta^t T(x) - \phi(\eta)},
\]
where
\[
\phi(\eta) = \log {\int_{\mathcal{X}} e^{\eta^t T(x)}
  d\nu(x)}
\]
is the logarithm of the Laplace transform of the measure $\nu^T=\nu\circ
T^{-1}$ that the statistic $T$ induces on the Borel $\sigma$-algebra of
$\RRR^k$.  The support of $\nu^T$ is the intersection of all closed sets
$A\subseteq \RRR^k$ that satisfy $\nu^T(\RRR^k\setminus A)=0$.  Recall that
the affine dimension of $A\subseteq \RRR^k$ is the dimension of the linear
space spanned by all differences $x-y$ of two vectors $x,y\in A$.

\begin{definition}
  Let $P_\eta$ be the probability measure on $(\mathcal{X},\mathcal{A})$
  that has $\nu$-density $p_\eta$.  The probability distributions
  $(P_\eta\mid \eta\in N)$ form a regular exponential family of order $k$
  if $N$ is an open set in $\RRR^k$ and the affine dimension of the support
  of $\nu^T$ is equal to $k$.  The statistic $T(x)$ that induces the
  regular exponential family is called a canonical sufficient statistic.
\end{definition}

The order of a regular exponential family is unique and if the same family
is represented using two different canonical sufficient statistics then
those two statistics are non-singular affine transforms of each other
\citep[Thm.~1.9]{brown:1986}.
\\

\par

\noindent {\bf 2.1. Examples}

Regular exponential families comprise families of discrete distributions,
which were the subject of much of the work on algebraic statistics.

\begin{example}[Discrete data]
  \label{ex:discrete}
  \rm Let the sample space $\mathcal{X}$ be the set of integers
  $\{1,\dots,m\}$.  Let $\nu$ be the counting measure on $\mathcal{X}$,
  i.e., the measure $\nu(A)$ of $A\subseteq \mathcal{X}$ is equal to the
  cardinality of $A$.  Consider the statistic $T:\mathcal{X}\to
  \RRR^{m-1}$,
  \[
  T(x) = \big(I_{\{1\}}(x),\dots,I_{\{m-1\}}(x)\big)^t,
  \]
  whose zero-one components indicate which value in $\mathcal{X}$ the
  argument $x$ is equal to.  In particular, when $x = m$, $T(x)$ is the zero vector.  The induced measure $\nu^T$ is a measure on the Borel
  $\sigma$-algebra of $\RRR^{m-1}$ with support equal to the $m$ vectors in
  $\{0,1\}^{m-1}$ that have at most one non-zero component.  The
  differences of these $m$ vectors include all canonical basis vectors of
  $\RRR^{m-1}$.  Hence, the affine dimension of the support of $\nu^T$ is
  equal to $m-1$.
  
  It holds for all $\eta\in\RRR^{m-1}$ that
  \[
  \phi(\eta) = 
  \log\left(1+\sum_{x=1}^{m-1} e^{\eta_x}\right)<\infty
  \]
  Hence, the natural parameter space $N$ is equal to all of $\RRR^{m-1}$
  and in particular is open.  The $\nu$-density $p_\eta$ is a probability
  vector in $\RRR^m$.  The components $p_\eta(x)$ for $1\le x\le m-1$ are
  positive and given by
  \[
  p_\eta(x) = \frac{e^{\eta_x}}{1+\sum_{x=1}^{m-1} e^{\eta_x}}.
  \]
  The last component of $p_\eta$ is also positive and equals
  \[
  p_\eta(m)=1-\sum_{x=1}^{m-1}
  p_\eta(x)=\frac{1}{1+\sum_{x=1}^{m-1} e^{\eta_x}}.
  \]
  The family of induced probability distribution $(P_\eta\mid
  \eta\in\RRR^{m-1})$ is a regular exponential family of order $m-1$. The
  interpretation of the natural parameters $\eta_x$ is one of log odds
  because $p_\eta$ is equal to a given positive probability vector
  $(p_1,\dots,p_m)$ if and only if $\eta_x = \log (p_x/p_m)$ for $x =
  1,\dots,m-1$.  This establishes a correspondence between the natural
  parameter space $N=\RRR^{m-1}$ and the interior of the $m-1$ dimensional
  probability simplex. \qed
\end{example}

The other distributional framework that has seen application of algebraic
geometry is that of multivariate normal distributions.

\begin{example}[Normal distribution]
  \label{ex:normal}
  \rm Let the sample space $\mathcal{X}$ be Euclidean space $\RRR^p$
  equipped with its Borel $\sigma$-algebra and Lebesgue measure $\nu$.
  Consider the statistic $T:\mathcal{X}\to\RRR^p\times\RRR^{p(p+1)/2}$
  given by
  \[
  T(x) =
  (x_1,\dots,x_p,-x_1^2/2,\dots,-x_p^2/2,-x_1x_2,\dots,-x_{p-1}x_p)^t.
  \]
  The polynomial functions that form the components of $T(x)$ are linearly
  independent and thus the support of $\nu^T$ has the full affine dimension
  $p+p(p+1)/2$.
  
  If $\eta\in\RRR^p\times\RRR^{p(p+1)/2}$, then write $\eta_{[p]}\in\RRR^p$
  for the vector of the first $p$ components $\eta_i$, $1\le i\le p$.
  Similarly, write $\eta_{[p\times p]}$ for the symmetric $p\times
  p$-matrix formed from the last $p(p+1)/2$ components $\eta_{ij}$, $1\le
  i\le j\le p$.  The function $x\mapsto e^{\eta^tT(x)}$ is $\nu$-integrable
  if and only if $\eta_{[p\times p]}$ is positive definite.  Hence, the
  natural parameter space $N$ is equal to the Cartesian product of $\RRR^p$
  and the cone of positive definite $p\times p$-matrices.  If $\eta$ is in
  the open set $N$, then
  \[
  \phi(\eta) = -\frac{1}{2}\left(\log \det(\eta_{[p\times p]}) -
  \eta_{[p]}^t\eta_{[p\times p]} \eta_{[p]}-p\log(2\pi)\right).
  \]
  The Lebesgue densities $p_\eta$ can be written as
  \[
  p_\eta(x) = \frac{1}{\sqrt{(2\pi)^{p}\det(\eta_{[p\times
        p]}^{-1})}}\exp\left\{ \eta_{[p]}^t x - \trace(\eta_{[p\times
      p]}xx^t)/2-\eta_{[p]}^t\eta_{[p\times p]} \eta_{[p]}/2 \right\}.
  \]
  Setting $\Sigma=\eta_{[p\times p]}^{-1}$ and $\mu=\eta_{[p\times
    p]}^{-1}\eta_{[p]}$, we find that
  \[
  p_\eta(x) =
  \frac{1}{\sqrt{(2\pi)^{p}\det(\Sigma)}}\exp\left\{-\tfrac{1}{2}
    (x-\mu)^t\Sigma^{-1}(x-\mu) \right\}
  \]
  is the density of the multivariate normal distribution
  $\ND_p(\mu,\Sigma)$.  Hence, the family of all multivariate normal
  distributions on $\RRR^p$ with positive definite covariance matrix is a
  regular exponential family of order $p+p(p+1)/2$.  \qed
\end{example}

The structure of a regular exponential family remains essentially unchanged
when sampling independent and identically distributed observations.

\begin{example}[Samples]  \rm
  A sample $X_1,\dots,X_n$ from $P_\eta$ comprises independent random
  vectors, all distributed according to $P_\eta$.  Denote their joint
  distribution by $\otimes_{i=1}^n P_\eta$.  An important property of a
  regular exponential family $(P_\eta\mid \eta\in N)$ of order $k$ is that
  the induced family $(\otimes_{i=1}^n P_\eta\mid \eta\in N)$ is again a
  regular exponential family of order $k$ with canonical sufficient
  statistic $\sum_{i=1}^n T(x_i)$ and Laplace transform $n\phi(\eta)$.  For
  discrete data as discussed in Example \ref{ex:discrete},
  the canonical sufficient statistic is given by the vector of counts
  \[
  N_x = \sum_{i=1}^n I_{\{x\}}(x_i),\quad x=1,\dots,m-1.
  \]
  For the normal distribution in Example \ref{ex:normal}, the canonical
  sufficient statistic is in correspondence with the empirical mean vector
  $\bar X$ and the empirical covariance matrix $S$; compare
  (\ref{eq:empiricalXS}).  \qed
\end{example}
\mbox{ }

\noindent {\bf 2.2. Likelihood inference in regular exponential families}

Among the nice properties of regular exponential families is their behavior
in likelihood inference.  Suppose the random vector $X$ is distributed
according to some unknown distribution from a regular exponential family
$(P_\eta\mid \eta\in N)$ of order $k$ with canonical sufficient statistic
$T$.  Given an observation $X=x$, the {\em log-likelihood function\/} takes
the form
\[
\ell(\eta\mid T(x))= \eta^t T(x) - \phi(\eta).
\]
The log-Laplace transform $\phi$ is a strictly convex and smooth, that is,
infinitely many times differentiable, function on the convex set $N$
\citep[Thm.~1.13, Thm.~2.2, Cor.~2.3]{brown:1986}.  The derivatives of
$\phi$ yield the moments of the canonical sufficient statistic such as the
expectation and covariance matrix,
\begin{align}
\label{eq:meanmap}
\zeta(\eta) &:= \displaystyle\frac{d}{d\eta}\phi(\eta) = \E_\eta[T(X)],\\
\nonumber
\Sigma(\eta) &:= \displaystyle\frac{d^2}{d\eta^2}\phi(\eta) =
\E_\eta\left\{[T(X)-\zeta(\eta)]\,
[T(X)-\zeta(\eta)]^t\right\}. 
\end{align}
The matrix $\Sigma(\eta)$ is positive definite since the components of $T(X)$
may not exhibit a linear relationship that holds almost everywhere.

The strict convexity of $\phi$ implies strict concavity of the
log-likelihood function $\ell$.  Hence, if the {\em maximum likelihood
  estimator\/} (MLE)
\[
\hat\eta(T(x)) = \arg\max_{\eta\in N} \ell(\eta\mid T(x))
\]
exists then it is the unique local and global maximizer of $\ell$ and can be
obtained as the unique solution of the likelihood equations $\zeta(\eta) =
T(x)$.  The existence of $\hat\eta(T(x))$ is equivalent to the condition
$T(x)\in\zeta(N)$; the open set $\zeta(N)$ is equal to the interior of the
convex hull of the support of $\nu^T$ \citep[Thm.~5.5]{brown:1986}.

If $X_1,\dots,X_n$ are a sample of random vectors drawn from
$P_\eta$, then the previous discussion applies to the family
$(\otimes_{i=1}^n P_\eta\mid \eta\in N)$.  In particular, the likelihood
equations become
\[
n \zeta(\eta) = \sum_{i=1}^n T(X_i) \iff \zeta(\eta) =
\bar T := \frac{1}{n}\sum_{i=1}^n T(X_i).
\]
By the strong law of large numbers, $\bar T$ converges almost surely to the
true parameter point $\zeta(\eta_0)\in \zeta(N)$.  It follows that the
probability of existence of the MLE, $\Prob_{\eta_0}\big(\bar T\in
\zeta(N)\big)$, tends to one as the sample size $n$ tends to infinity.
Moreover, the {\em mean parametrization\/} map $\eta\mapsto \zeta(\eta)$ is
a bijection from $N$ to $\zeta(N)$ that has a differentiable inverse with
total derivative
\[
\frac{d}{d\eta}\zeta^{-1}(\eta) = \Sigma(\eta)^{-1},
\]
which implies in conjunction with an application of the central limit
theorem:
\begin{proposition}
  The MLE $\hat\eta(\bar T)=\zeta^{-1}(\bar T)$ in a regular exponential
  family is asymptotically normal in the sense that if $\eta_0$ is the true
  parameter, then 
  \[
  \sqrt{n}[\hat\eta(\bar T) - \eta_0]
  \stackrel{n\to\infty}{\longrightarrow_d} \ND_k (0,\Sigma(\eta_0)^{-1}). 
  \]
\end{proposition}


A submodel of a regular exponential family $(P_\eta\mid\eta\in N)$ of order
$k$ is given by a subset $M\subseteq N$.  If the geometry of the set $M$ is
regular enough, then the submodel may inherit the favorable properties of
likelihood inference from its reference model, the regular exponential
family.  The nicest possible case occurs when the submodel
$(P_\eta\mid\eta\in M)$ has parameter space $M=N\cap L$, where $L\subseteq
\RRR^k$ is an affine subspace of $\RRR^k$.  Altering the canonical
sufficient statistics one finds that $(P_\eta\mid\eta\in M)$ forms a
regular exponential family of order $\dim(L)$.  

Given a single observation $X$ from $P_{\eta}$, the {\em likelihood ratio
  test\/} for testing $H_0:\eta\in M$ versus $H_1:\eta\in N\setminus M$
rejects $H_0$ for large values of the {\em likelihood ratio statistic\/}
\[
\lambda_M(T(X)) = \sup_{\eta\in N} \ell(\eta\mid T(X)) - \sup_{\eta\in M}
\ell(\eta\mid T(X)).
\]
If we observe a sample $X_1,\dots,X_n$ from $P_{\eta}$, then the likelihood
ratio statistic depends on $\bar T$ only and is equal to $n\lambda_M(\bar
T)$.  For a rejection decision, the distribution of $n\lambda_M(\bar T)$
can often be approximated using the next asymptotic result.

\begin{proposition}
  If $M=N\cap L$ for an affine space $L$ and the true parameter $\eta_0$ is
  in $M$, then the likelihood ratio statistic $n\lambda_M(\bar T)$
  converges to $\chi^2_{k-\dim(L)}$, the chi-square distribution with
  $k-\dim(L)$ degrees of freedom, as $n\to\infty$.
\end{proposition}

In order to obtain asymptotic results such as uniformly valid chi-square
asymptotics for the likelihood ratio statistic, the set $M$ need not be
given by an affine subspace.  In fact, if $M$ is an $m$-dimensional smooth
manifold in $\RRR^k$, then $n\lambda_M(\bar T)$ still converges in
distribution to $\chi^2_{k-m}$ for any $\eta_0\in M$.  A set $M$ is an
$m$-dimensional {\em smooth manifold\/} if for all $\eta_0\in M$ there
exists an open set $U\subseteq \RRR^k$ containing $\eta_0$, an open set
$V\subseteq \RRR^k$, and a diffeomorphism $g:V\to U$ such that
$g\big(V\cap(\RRR^m\times\{0\})\big) = U$.  Here,
$\RRR^m\times\{0\}\subseteq \RRR^k$ is the subset of vectors for which the
last $k-m$ components are equal to zero.  A {\em diffeomorphism\/} $g:V\to
U$ is a smooth bijective map that has a smooth inverse $g^{-1}:U\to V$.  An
exponential family induced by a smooth manifold in the natural parameter
space is commonly termed a {\em curved exponential family\/}; see
\cite{kass:1997} for an introduction to this topic.

The fact that many interesting statistical models, in particular models
involving hidden variables, are not curved exponential families calls for
generalization.  One attempt at such generalization was made by
\cite{geiger:2001} who introduce so-called {\em stratified exponential
  families\/}.  A stratified exponential family is obtained by piecing
together several curved exponential families.  However, as the next example
shows, stratified exponential families appear to be a bit too general as a
framework unless more conditions are imposed on how the curved exponential
families are joined together.  Example \ref{ex:stratified} is inspired by
an example in \citet[p.~199]{rockafellar}.

\begin{example}
  \label{ex:stratified}
  \rm Consider the regular exponential family $\mathcal{P}$ of bivariate
  normal distributions with unknown mean vector $\mu=(\mu_1,\mu_2)^t\in
  \RRR^2$ but covariance matrix $\Sigma$ equal to the identity matrix
  $I_2\in\mathbb{R}^{2\times 2}$.  The natural parameter space of this
  model is the plane $\RRR^2$.  When drawing a sample $X_1,\dots,X_n$ from
  a distribution in $\mathcal{P}$, the canonical statistic is the sum of
  the random vectors.  Dividing by the sample size $n$ yields the sample
  mean vector $\bar X\in\RRR^2$, which is also the MLE of $\mu$.  In the
  following we will assume that the true parameter $\mu_0$ is equal to the
  origin.  Then the rescaled sample mean vector $\sqrt{n}\,\bar X$ is
  distributed according to the bivariate standard normal distribution
  $\mathcal{N}_2(0,I_2)$.
  
  If we define a submodel $\mathcal{P}_C\subseteq\mathcal{P}$ by
  restricting the mean vector to lie in a closed set $C\subseteq \RRR^2$,
  then the MLE $\hat\mu$ for the model $\mathcal{P}_C$ is the point in $C$
  that is closest to $\bar X$ in Euclidean distance.  For $n=1$, the
  likelihood ratio statistic $\lambda_C(\bar X)$ for testing $\mu\in
  C$ versus $\mu\not\in C$ is equal to the squared Euclidean distance
  between $\bar X$ and $C$.  Hence, the likelihood ratio statistic based on
  an $n$-sample is
  \[
  n\lambda_C(\bar X) = n\cdot \min_{\mu\in C} || \bar X - \mu ||^2
  = \min_{\mu\in \sqrt{n}\, C} || \sqrt{n}\,\bar X - \mu ||^2,
  \]
  i.e., the squared Euclidean distance between the standard normal random
  vector $\sqrt{n}\,\bar X$ and the rescaled set $\sqrt{n}\,C$.
 
  As a concrete choice of a submodel, consider the set
  \[
  C_1 = \{ (\mu_1,\mu_2)^t\in\RRR^2\mid \mu_2=\mu_1\sin(1/\mu_1),\;
  \mu_1\not=0\}\cup \{(0,0)^t\}.
  \]
  This set is the disjoint union of the two one-dimensional smooth
  manifolds obtained by taking $\mu_1<0$ and $\mu_1>0$, and the
  zero-dimensional smooth manifold given by the origin.  These manifolds
  form a stratification of $C_1$ \citep[p.~513]{geiger:2001}, and thus the
  model $\mathcal{P}_{C_1}$ constitutes a stratified exponential family.
  In Figure \ref{fig:tinvgraph}, we plot three of the sets $\sqrt{n}\,C_1$
  for the choices $n=100,100^2,100^3$.  The range of the plot is restricted
  to the square $[-3,3]^2$, which contains the majority of the mass of the
  bivariate standard normal distribution.  The figure illustrates the fact
  that as $n$ tends to infinity the sets $\sqrt{n}\,C_1$ fill more and more
  densely the 2-dimensional cone comprised between the axes
  $\mu_2=\pm\mu_1$.  Hence, $n\lambda_{C_1}(\bar X)$ converges in
  distribution to the squared Euclidean distance between a bivariate
  standard normal point and this cone.  So although we pieced
  together smooth manifolds of codimension 1 or larger, the limiting
  distribution of the likelihood ratio statistic is obtained from a
  distance to a full-dimensional cone.

  \begin{figure}[t]
    \centering
    \begin{tabular}{lll}
      \includegraphics[width=4cm]{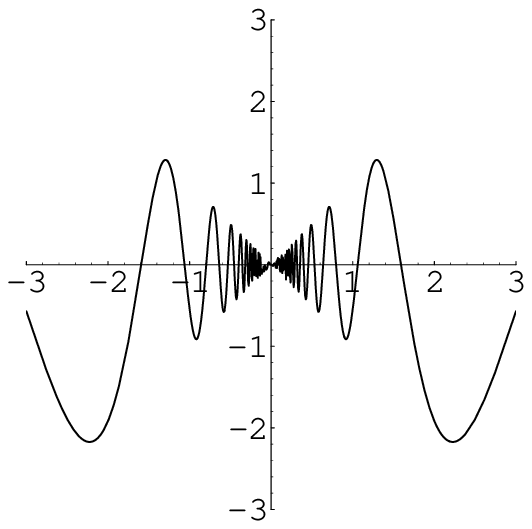}  &
      \includegraphics[width=4cm]{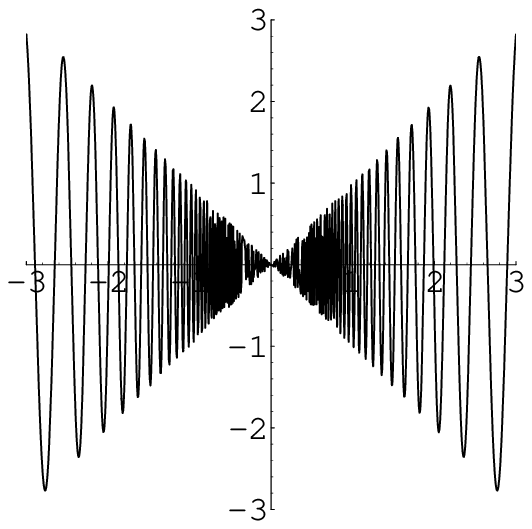}  &
      \includegraphics[width=4cm]{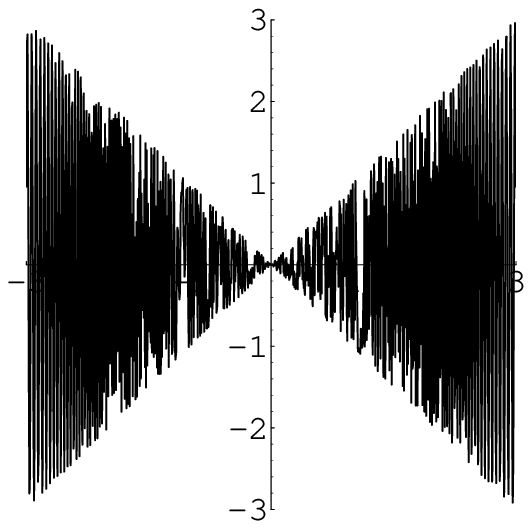}  
    \end{tabular}
    \caption{Sets $\sqrt{n}\,C_1$ for $n=100,100^2,100^3$.}
    \label{fig:tinvgraph}
  \end{figure}
  
  As a second submodel consider the one induced by the set
  \[
  C_2 = \{ (\mu_1,\mu_2)^t\in\RRR^2\mid \mu_2=\mu_1\sin(-\log(|\mu_1|/4)),\;
  \mu_1\in [-3,3]\setminus\{0\}\}\cup \{(0,0)^t\}.
  \]
  The model $\mathcal{P}_{C_2}$ is again a stratified exponential family.
  However, now the sets $\sqrt{n}\,C_2$ have a wave-like structure even for
  large sample sizes $n$; compare Figure \ref{fig:logtgraph}.  We conclude
  that in this example the likelihood ratio test statistic
  $n\lambda_{C_2}(\bar X)$ does not converge in distribution as $n$ tends to
  infinity. \qed
  \begin{figure}[t]
    \centering
    \begin{tabular}{lll}
      \includegraphics[width=4cm]{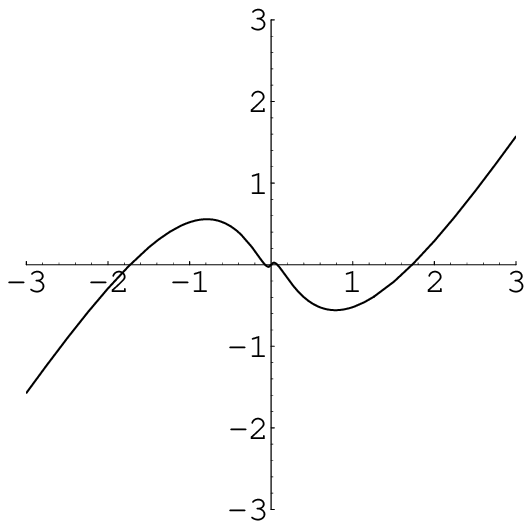}  &
      \includegraphics[width=4cm]{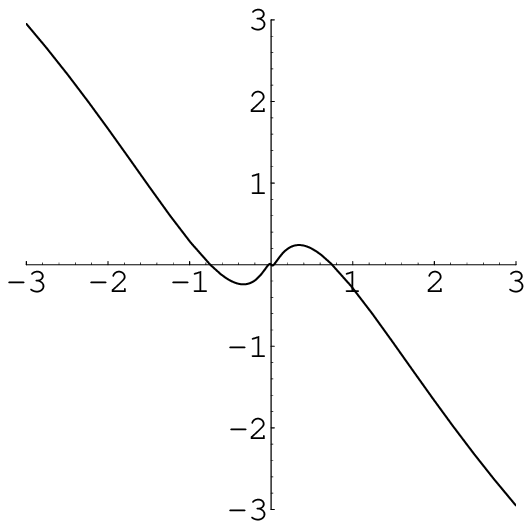}  &
      \includegraphics[width=4cm]{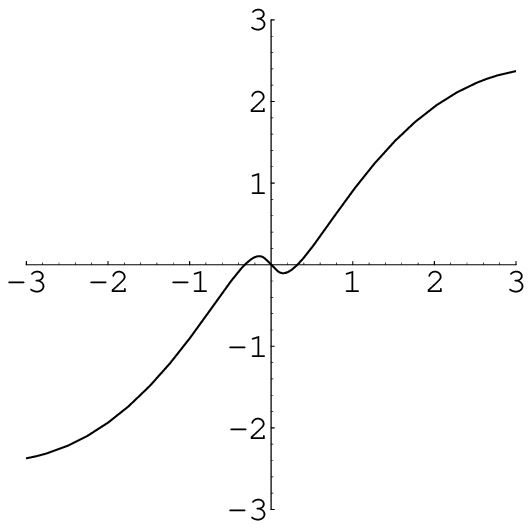}  
    \end{tabular}
    \caption{Sets $\sqrt{n}\,C_2$ for $n=100,100^2,100^3$.}
    \label{fig:logtgraph}
  \end{figure}
\end{example}

The failure in the previous example of nice asymptotic behavior of the
likelihood ratio test is part of our motivation for restricting to the
class of \emph{algebraic exponential families}, which we introduce 
next.
\\

\par


\setcounter{section}{3}
\setcounter{equation}{0} 

\noindent {\bf 3. Algebraic exponential families}

In the following definition, which was anticipated in the introduction, we
propose the use of semi-algebraic sets to unify different definitions of
algebraic statistical models.  Using semi-algebraic sets eliminates
phenomena as created in Example \ref{ex:stratified} because these sets have
nice local geometric properties.  In addition, imposing algebraic structure
allows one to employ the tools of computational algebraic geometry to
address questions arising in statistical inference.  (More details on both
these points are given in Section 4.)

\newcounter{olddefinition}
\setcounter{olddefinition}{\value{definition}}
\setcounter{definition}{1}
\begin{definition}
Let $(P_\eta\mid\eta\in N)$ be a regular exponential family of order $k$.
  The subfamily induced by the set $M\subseteq N$ is an algebraic
  exponential family if there exists an open set $\bar N\subseteq \RRR^k$,
  a diffeomorphism $g:N\to \bar N$, and a semi-algebraic set $A \subseteq
  \RRR^k$ such that $M=g^{-1}( A \cap \bar N)$.
\end{definition}
\setcounter{definition}{\value{olddefinition}}


The definition states that an algebraic exponential family is given by a
semi-algebraic subset of the parameter space of a regular exponential
family.  However, this parameter space may be obtained by a
reparametrization $g$ of the natural parameter space $N$, which provides
the necessary flexibility to capture the algebraic structure found in
interesting statistical models including ones that do not form curved
exponential families.  The mean parametrization $\zeta(\eta)$ is one
example of a useful reparametrization.  

Before giving examples of algebraic exponential families we provide some
background on semi-algebraic sets; more in depth introductions can be
found, for example, in \cite{benedetti:1990} or \cite{bochnak:1998}.
\\

\par


\noindent {\bf 3.1. Basic facts about semi-algebraic sets}

A monomial in indeterminates (polynomial variables) $t_1, \ldots, t_n$, is
a formal expression of the form $\bft^\beta = t_1^{\beta_1} t_2^{\beta_2}
\cdots t_n^{\beta_n}$ where $\beta = (\beta_1, \ldots, \beta_n)$ is the
non-negative integer vector of exponents.  A polynomial
$$f = \sum_{ \beta \in B} c_\beta \bft^\beta$$
is a linear combination of
monomials where the coefficients $c_\beta$ are in a fixed field $\bbk$ and
$B \subset \bbn^n$ is a finite set of exponent vectors.  The collection of
all polynomials in the indeterminates $t_1, \ldots, t_n$ with coefficients
in a fixed field $\bbk$ is the set $\bbk[\bft] = \bbk[t_1, \ldots, t_n]$.
The collection of polynomials $\bbk[\bft]$ has the algebraic structure of a
\emph{ring}.  Each polynomial in $\bbk[\bft]$ is a formal linear
combination of monomials that can also be considered as a function $f:
\bbk^n \rightarrow \bbk$, defined by evaluation.  Throughout the paper, we
will focus attention on the ring $\bbr[\bft]$ of polynomials with real
coefficients.

\begin{definition}
A \emph{basic semi-algebraic set} is a subset of points in $\bbr^n$ of the form
$$
A = \{ \theta \in \bbr^n \mid f(\theta) > 0\; \forall f \in
F,\; h(\theta) = 0\; \forall h \in H \}
$$
where $F \subset \bbr[\bft]$ is a finite (possibly empty) collection of
polynomials and $H\subseteq \bbr[\bft]$ is an arbitrary (possibly empty)
collection of polynomials.  A \emph{semi-algebraic set} is a finite union
of basic semi-algebraic sets.  If $F = \emptyset$ then $A$ is called a
\emph{real algebraic variety}.
\end{definition}

A particular special case of a general 
semi-algebraic set occurs when we consider sets of the form
$$
A = \{ \theta \in \bbr^n \mid f(\theta) > 0\; \forall f \in F,\; g(\theta)
\geq 0\; \forall g \in G,\; h(\theta) = 0\; \forall h \in H \}
$$
where both $F$ and
$G$ are finite collections of real polynomials.  

\begin{example}  \rm
  The open probability simplex for discrete random variables is a basic
  semi-algebraic set, where $F = \{t_i \mid i = 1, \ldots, n-1\}\cup\{ 1-
  \sum_{i=1}^{n-1} t_i \}$ and $H = \emptyset $.  More generally, the
  relative interior of any convex polyhedron in any dimension is a basic
  semi-algebraic set, while the whole polyhedron is an ordinary
  semi-algebraic set.  \qed
\end{example}

\begin{example}  \rm
  The set $\Sigma \subset \bbr^{m \times m}$ of positive definite matrices
  is a basic semi-algebraic set, where $F$ consists of all principal
  subdeterminants of a symmetric matrix $\Psi$, and $G$ is the empty set.
 \qed
\end{example} 

In our introduction, parametrically specified statistical models were
claimed to be algebraic statistical models.  This non-trivial claim holds
due to the famous Tarski-Seidenberg theorem \citep[e.g.][]{bochnak:1998},
which says that the image of a semi-algebraic set under any nice enough
mapping is again a semi-algebraic set.  To make this precise we need to
define the class of mappings of interest.


Let $\psi_1 = f_1/g_1, \ldots, \psi_n = f_n/g_n$ be rational functions
where $f_i,g_i \in \bbr[\bft] = \bbr[t_1, \ldots, t_d]$ are real polynomial
functions.  These rational functions can be used to define a rational map
$$\psi: \bbr^d \rightarrow \bbr^n, \quad \bfa \mapsto (\psi_1(\bfa),
\ldots, \psi_n(\bfa)),$$
which is
well-defined
on the open set $D_\psi =\{ \bfa \subset \bbr^d : \prod g_i(\bfa) \neq 0
\}$.

\begin{theorem}[Tarski-Seidenberg]
  Let $A \subseteq \bbr^d$ be a semi-algebraic set and $\psi$ a
  rational map that is well-defined on $A$, that is,
  $A  \subseteq D_\psi$.  Then the image $\psi(A)$ is a
  semi-algebraic set.
\end{theorem}

\cite{ascb} define an algebraic statistical model as the image of a
polynomial parametrization $\psi(A) \subseteq \Delta$ where $A$ is the
interior of a polyhedron and $\Delta$ is the probability simplex.  The
emphasis on such models, which one might call \emph{parametric algebraic
  statistical models}, results from the fact that most models used in the
biological applications under consideration (sequence alignment and
phylogenetic tree reconstruction, to name two) are parametric models for
discrete random variables.  Furthermore, the precise algebraic form of
these parametric models is essential to parametric maximum a posteriori
estimation, one of the major themes in the text of \cite{ascb}.  The
Tarski-Seidenberg theorem and Example \ref{ex:discrete} yield the
following unifying fact.


\begin{corollary}
  If a parametric statistical model for discrete random variables is a
  well-defined image of a rational map from a semi-algebraic set to the
  probability simplex, then the model is an algebraic exponential family.
\end{corollary}


\par
\noindent {\bf 3.2.  Independence models as examples}

Many statistical models are defined based on considerations of
(conditional) independence.  Examples include Markov chain models, models
for testing independence hypotheses in contingency tables and graphical
models, see e.g.~\cite{lauritzen:1996}.  As we show next, conditional
independence yields algebraic exponential families in both the Gaussian and
discrete cases.  The algebraic structure also passes through under
marginalization, as we will illustrate in Section 4.

\begin{example}[Conditional independence in normal distributions]
  \label{ex:cigauss}
  \rm Let $X=(X_1,\dots,X_p)$ be a random vector with joint normal
  distribution $\mathcal{N}_p(\mu,\Sigma)$ with mean vector $\mu\in\RRR^p$
  and positive definite covariance matrix $\Sigma$.  For three pairwise
  disjoint index sets $A,B,C\subseteq \{1,\dots,p\}$, the subvectors $X_A$
  and $X_B$ are conditionally independent given $X_C$, in symbols $X_A\ind
  X_B\mid X_C$ if and only if
  \begin{equation*}
    \det(\Sigma_{\{i\}\cup C\times
      \{j\} \cup C}) = 0 \quad\forall i\in A,\; j\in B.
  \end{equation*}
  If $C=\emptyset$, then conditional independence given $X_\emptyset$ is
  understood to mean marginal independence of $X_A$ and $X_B$.  
   \qed
\end{example}

\begin{example}[Conditional independence in the discrete case] \rm
\label{ex:cidisc}
Conditional independence statements also have a natural algebraic
interpretation in the discrete case.  As the simplest example, consider the
conditional independence statement $X_1 \ind X_2 \mid X_3$ for the discrete
random vector $(X_1, X_2, X_3)$.  This translates into the collection of
algebraic constraints on the joint probability distribution
\begin{eqnarray*}
&  &  \Prob(X_1 = i_1, X_2 = j_1, X_3 = k) \cdot 
      \Prob(X_1 = i_2, X_2 = j_2, X_3 = k)   \\
 & =  & \Prob(X_1 = i_1, X_2 = j_2, X_3 = k) \cdot 
        \Prob(X_1 = i_2, X_2 = j_1, X_3 = k)  
\end{eqnarray*} 
for all $i_1, i_2 \in [m_1]$, $j_1, j_2 \in [m_2]$ and $k \in [m_3]$, where
$[m] = \{1,2,\ldots, m\}$.  Alternatively, we might write this in a more
compact algebraic way as:
$$p_{i_1j_1k} p_{i_2j_2k} - p_{i_1j_2k} p_{i_2j_1k} \, \, = \, \, 0,$$
where $p_{ijk}$ is shorthand for $\Prob(X_1 = i, X_2 = j, X_3 = k)$.  In
general, any collection of conditional independence statements for discrete
random variables corresponds to a collection of quadratic polynomial
constraints on the components of the joint probability vector.  \qed
\end{example}  
\mbox{ } 
\par

\setcounter{section}{4}
\setcounter{equation}{0} 


\noindent {\bf 4. Model geometry} 

Of fundamental importance to statistical inference is the intuitive notion
of the ``shape'' of a statistical model, reflected in its abstract
geometrical properties.  Examples of interesting geometrical features are
whether or not the likelihood function is multimodal, whether or not the
model has singularities (is non-regular) and the nature of the underlying
singularities.  These are all part of answering the question: How does the
geometry of the model reflect its statistical features?  When
the model is an algebraic exponential family, these problems can be
addressed using algebraic techniques, in particular by computing with
ideals.  This is even true when the model comes in a parametric form,
however, it is then often helpful to translate to an
implicit representation of the model.\\

\par


\noindent {\bf 4.1 Model invariants}

Recall that an ideal $I \subset \bbr[\bft]$ is a collection of polynomials
such that for all $f, g \in I$, $f + g \in I$ and for all $f \in I$ and $h
\in \bbr[\bft]$, $h \cdot f \in I$.  Ideals can be used to determine real
algebraic varieties by computing the zero set of the ideal:
$$V(I) \quad = \quad \left\{ \bfa \in \bbr^n \, \, | \, \, f(\bfa) = 0
  \mbox{ for all } f \in I \right\}.$$
When we wish to speak of the variety
over the complex numbers we use the notation $V_\bbc(I)$.  Reversing this
procedure, if we are given a set $V \subset \bbr^n$ we can compute its
defining ideal, which is the set of all polynomials that vanish on $V$:
$$I(V) \quad = \quad \left\{ f \in \bbr[\bft] \, \, | \, \, f(\bfa) = 0
  \mbox{ for all } \bfa \in V \right\}.$$

\begin{definition}
  Let $A$ be a semi-algebraic set defining an algebraic exponential family
  $\mathcal{P}_M=(P_\eta\mid\eta\in M)$ via $M=g^{-1}(A\cap g(N))$.  A
  polynomial $f$ in the vanishing ideal $I(A)$ is a {\em model invariant}
  for $\mathcal{P}_M$.
\end{definition}

\begin{remark}
  \rm The term ``model invariant'' is chosen in analogy to the term
  ``phylogenetic invariant'' that was coined by biologists working with
  statistical models that are useful for the reconstruction of phylogenetic
  trees.
\end{remark}

Given a list of polynomial $f_1, \ldots, f_k$ the ideal generated by these
polynomials is denoted
$$\left< f_1, \ldots, f_k \right> \quad = \quad \left\{ \sum_{i =1}^k h_i
  \cdot f_i \, \, | \, \, h_i \in \bbr[\bft] \right\}.$$
The Hilbert basis
theorem says that every ideal in a polynomial ring has a finite generating
set.  Thus, when working with a statistical model that we want to describe
algebraically, we need to compute a \emph{finite list} of polynomials that
generate the ideal of model invariants.  These equations can be used to
address questions like determining the structure of singularities which in
turn can be used to address asymptotic questions.

\begin{example} [Conditional independence]\rm \label{ex:ciideal}
  In Example \ref{ex:cigauss} we gave a set of equations whose zero set in
  the cone of positive definite matrices is the independence model obtained
  from $X_A\ind X_B\mid X_C$.  However, there are more equations, in
  general, that belong to the ideal of model invariants $I$.  In
  particular, we have
  $$I \quad = \quad \left< \det \tilde\Sigma \, \, | \,\, \tilde\Sigma
    \mbox{ is a } (|C| + 1) \times (|C| + 1) \mbox{ submatrix of }
    \Sigma_{A \cup C, B \cup C} \right>.$$
  The fact that this ideal
  vanishes on the model follows from the fact that any $\Sigma$ in the
  model is positive definite and, hence, each principal minor is
  invertible.  The fact that the indicated ideal comprises all model
  invariants can be derived from a result in commutative algebra
  \citep{Conca1994}.
  
  In the discrete case, the polynomials we introduced in Example
  \ref{ex:cidisc} generate the ideal of model invariants for the model
  induced by $X_1\ind X_2\mid X_3$.  For models induced by collections of
  independence statements this need no longer be true; compare Theorem 8 in
  \cite{garcia:2005}. \qed
\end{example}

One may wonder what the use of passing from the set of polynomials
exhibited in Example \ref{ex:cigauss} to the considerably larger set of
polynomials described in Example \ref{ex:ciideal} is, since both sets of
polynomials define the model inside the cone of positive definite matrices.
The smaller set of polynomials have the property that there lie singular
covariance matrices in the positive \emph{semi}definite cone that satisfy
the polynomial constraints but are not limits of covariance matrices in the
model.  From an algebraic standpoint, the main problem is that the ideal
generated by the smaller set of polynomials is not a prime ideal.  In
general, we prefer to work with the prime ideal given by all model
invariants because prime ideals tend to be better behaved from a
computational standpoint and are less likely to introduce extraneous
solutions on boundaries.

For the conditional independence models described thus far, the equations
$I(A)$ that define the model come from the definition of the model.  For
instance, conditional independence imposes natural constraints on
covariance matrices of normal random variables and the joint probability
distributions of discrete random variables.  When we are presented with a
parametric model, however, it is in general a challenging problem of
computational algebra to compute the \emph{implicit} description of the
model $A$ as a semi-algebraic set.  At the heart of this problem is the
computation of the ideal of model invariants $I(A)$, which can be solved
using Gr\"{o}bner bases.  Methods for computing an implicit description
from a parametric description can be found in \cite{Cox1997}, though the
quest for better implicitization methods is an active area of research.

The vanishing ideal of a semi-algebraic set can be used to address many
questions about it, for instance, the dimension of a semi-algebraic set.
The following definition and proposition provide a useful characterization
of the dimension of a semi-algebraic set.


\begin{definition}
  A set of indeterminates $p_{i_1}, \ldots, p_{i_k}$ is \emph{algebraically
    independent} for the ideal $I$ if there is no polynomial only in
  $p_{i_1}, \ldots, p_{i_k}$ that belongs to $I$.
\end{definition}

\begin{proposition}\label{prop:dim}
  The dimension of a semi-algebraic set $A$ is the cardinality of the
  largest set of algebraically independent indeterminates for $I(A)$.
\end{proposition}

The proof that algebraically independent sets of indeterminates and Proposition
\ref{prop:dim} meshes with the usual geometric notion of dimension can be
found in \cite{Cox1997}. 

The subset $V_{\rm sing} \subset V$ where a variety $V$ is singular is also
a variety.  Indeed, suppose that polynomials $f_1, \ldots, f_k$ generate
the vanishing ideal $I(V)$.  Let $J \in \bbr[\bfx]^{k \times n}$ denote the
Jacobian matrix with entry $J_{ij} = \frac{\partial f_i}{\partial x_j}$.

\begin{proposition} \label{prop:sing}
  A point $\bfa \in V_\bbc(I)$ is a singular point of the complex
  variety if and only if $J(\bfa)$ has rank less than the codimension of the
  largest irreducible component of $V$ containing $\bfa$.
\end{proposition}

The singularities of the real variety are defined to be the intersection of
the singular locus of $V_\bbc(I)$ with $\bbr^n$.  Proposition
\ref{prop:sing} yields a direct way to compute, as an algebraic variety,
the singular locus of $V$.  Indeed, the rank of the Jacobian matrix is less
than $c$, if and only if the $c \times c$ minors of $J$ are all zero.
Thus, if $I$ defines an irreducible variety of codimension $c$, the ideal
$\left< M_c(J), f_1, \ldots, f_k \right>$ has as zero set the singular
locus of $V$, where $M_c(J)$ denotes the set of $c \times c$ minors of $J$.
If the variety is not irreducible, the singular set consists of the union
of the singular set of all the irreducible components together with the
sets of all pairwise intersections between irreducible components.

Removing the singularities $V_{\rm sing}$ from $V$ one obtains a smooth
manifold such that the local geometry at a non-singular point of $V$ is
determined by a linear space, namely, the tangent space.  At singular
points, the local geometry can be described using the
tangent cone, which is the semi-algebraic set that approximates the limiting
behavior of the secant lines that pass through the point of interest.  In
the context of parameter spaces of statistical models, the study of this
limiting behavior is crucial for the study of large sample asymptotics at a
singular point.  The geometry of the tangent cone for semi-algebraic sets
can be complicated and we postpone an in-depth study for a later
publication.  For the singular models that we encounter in the next
section, the crucial point on the tangent cone is the following
proposition.

\begin{proposition}
  \label{prop:tangentunion}
  Suppose that $A = V_1 \cup \cdots \cup V_m$ is the union of smooth
  algebraic varieties and let $\bfa$ be a point in the intersection $V_1
  \cap \cdots \cap V_j$ such that $\bfa\not\in V_k$ for $k\ge j+1$.  Then
  the tangent cone of $A$ at $\bfa$ is the union of the tangent planes to
  $V_1, \ldots, V_j$ at $\bfa$.
\end{proposition}

\par


\noindent {\bf 4.2 A conditional independence model with singularities}

Let $X=(X_1,X_2,X_3)$ have a trivariate normal distribution
$\ND_3(\mu,\Sigma)$, and define a model by requiring that $X_1\ind X_2$ and
simultaneously $X_1\ind X_2\mid X_3$.  By Example \ref{ex:cigauss}, the
model is an algebraic exponential family given by the subset
$M=\zeta^{-1}(A\cap \zeta(N))$, where
$\zeta(N)=\mathbb{R}^3\times\mathbb{R}^{3\times 3}_{\rm pd}$ is the
Gaussian mean parameter space and the algebraic variety
\[
A = \left\{ (\mu,\Sigma)\in\mathbb{R}^3\times\mathbb{R}^{3\times 3}_{\rm
    sym} \mid \sigma_{12}=0,\; \det(\Sigma_{\{1,3\}\times \{2,3\}})=
  \sigma_{12}\sigma_{33}-\sigma_{13}\sigma_{23}= 0 \right\}.
\]
Here, $\mathbb{R}^{3\times 3}_{\rm sym}$ is the space of symmetric $3\times
3$-matrices.  The set $A$ is defined equivalently by the joint vanishing of
$\sigma_{12}$ and $\sigma_{13}\sigma_{23}$. Hence, $A = A_{13}\cup
A_{23}$ for
\[
\begin{split}
A_{13} & =
\{(\mu,\Sigma)\in A\mid \sigma_{12}=\sigma_{13}=0\},\\
A_{23} &= 
\{(\mu,\Sigma)\in A\mid \sigma_{12}=\sigma_{23}=0\}.
\end{split}
\] 
This decomposition as a union reflects the well-known fact that
\[
\left[\;X_1\ind X_2\, \wedge\, X_1\ind X_2\mid X_3 \;\right] \iff
\left[\;X_1\ind (X_2,X_3)\, \vee\, X_2\ind (X_1,X_3) \;\right],
\]
which holds for the multivariate normal distribution but also when $X_3$ is
a binary variable; compare \cite[Thm.~8.3]{dawid:1980}.  By Proposition
\ref{prop:tangentunion} the singular locus of $A$ is the intersection
\[
A_{\rm sing} = A_{13}\cap A_{23} = \{(\mu,\Sigma)\in A\mid
\sigma_{12}=\sigma_{13}=\sigma_{23}=0\},
\]
which corresponds to diagonal covariance matrices $\Sigma$, or in other
words, complete independence of the three random variables $X_1\ind X_2\ind
X_3$.

Given $n$ independent and identically distributed normal random vectors
$X_1,\dots,X_n\in \RRR^3$, define the empirical mean and
covariance matrix as
\begin{equation}
\label{eq:empiricalXS}
\bar X = \frac{1}{n}\sum_{i=1}^n X_i, \quad
S = \frac{1}{n}\sum_{i=1}^n (X_i-\bar X)(X_i-\bar X)^t,
\end{equation}
respectively.  The likelihood ratio test statistic for testing the model
based on parameter space $M$ against the regular exponential family of all
trivariate normal distributions can be expressed as
\begin{equation}
  \label{eq:cilr2}
  \lambda_M(\bar X,S) = \log\left(
    \frac{s_{11}s_{22}}{s_{11}s_{22}-s_{12}^2}
  \right)
  +
  \min\left\{ \log\left(
      \frac{s_{33.2}}{s_{33.12}}\right),
    \log\left(
      \frac{s_{33.1}}{s_{33.12}}
    \right)\right\},
\end{equation}
where for $A\subseteq\{1,2\}$, $s_{33.A}$ is the empirical conditional
variance
\[
s_{33.A} = s_{33}-S_{\{3\}\times
  A}S_{A\times A}^{-1}S_{A\times \{3\}}.
\]
The three terms in
(\ref{eq:cilr2}) correspond to tests of the hypotheses
\[
X_1\ind X_2, \quad X_1\ind X_3\mid X_2, \quad\text{and}\quad X_2\ind
  X_3\mid X_1.
\]
Note that a joint distribution satisfies $X_1\ind (X_2,X_3)$ if and only if
it satisfies both $X_1\ind X_2$ and $X_1\ind X_3\mid X_2$.

If $(\mu,\Sigma)$ is an element of the smooth manifold $A\setminus A_{\rm
  sing}$, then $\lambda_M(\bar X,S)$ converges to a $\chi^2_2$-distribution
as $n$ tends to infinity; but over the singular locus the limiting
distribution is non-standard as detailed in \cite{drton:prag}.

\begin{proposition}
  \label{prop:cisingasy}
  Let $(\mu,\Sigma)\in A_{\rm sing}$.  As $n\to\infty$, the likelihood
  ratio test statistic $\lambda_M(\bar X,S)$ converges to the minimum of
  two dependent $\chi^2_2$-distributed random variables, namely,
  \[
  \lambda_M(\bar X,S) \longrightarrow_d
  \min(W_{12}+W_{13},W_{12}+W_{23})=
  W_{12}+\min(W_{13},W_{23})
  \]
  for three independent $\chi^2_1$-random variables $W_{12}$, $W_{13}$ and $W_{23}$.
\end{proposition}


Similar asymptotics arise in the model of joint marginal and conditional
independence in the discrete case with $X_3$ binary.  In this case the
variety breaks again into the union of two independence varieties $X_1 \ind
\{X_2, X_3 \}$ and $ X_2 \ind \{X_1, X_3 \}$, whose intersection is the
complete independence variety corresponding to $X_1 \ind X_2 \ind X_3$.
Non-standard asymptotics will occur at the intersection of these two
varieties.  However, as both of the varieties $X_1 \ind \{X_2, X_3\}$ and
$X_2 \ind \{X_1, X_3\}$ are smooth, the tangent cone is simply the union of
the two tangent spaces to the two component varieties.  The asymptotics
behave in a manner similar to the Gaussian case, as the minimum of chi-square
distributions.\\ 

\par


\noindent {\bf 4.3 Hidden random variables}

Another important use for the implicit equations defining a model are that
they can be used to determine a (partial) description of any new models
that arise from the given model via marginalization.  In particular,
algebraic methods can be used to explore properties of models with hidden
random variables.  In this section, we describe how to derive model
invariants via elimination in the presence of hidden variables for Gaussian
and discrete random variables.

\begin{proposition}
  Suppose that the random vector $X = (X_1, \ldots, X_p)$ is distributed
  according to a multivariate normal distribution from a model with ideal
  of model invariants $I \subset \bbr[\mu_i,\sigma_{ij} \, \, | \,\, 1 \leq
  i \leq j \leq p ]$.  Then the {\em elimination} ideal $I \cap
  \bbr[\mu_i,\sigma_{ij} \, \, | \, \, 1 \leq i \leq j \leq p-1]$ comprises
  the model invariants of the model created by marginalizing to $X' = (X_1,
  \ldots, X_{p-1})$.
\end{proposition}


The indicated elimination can be computed using Gr\"obner bases
\citep{Cox1997}.  A similar type of elimination formulation can be given
for the marginalization in the discrete case.

\begin{proposition}
  Let $X_1,\dots,X_p$ be discrete random variables with $X_k$ taking values
  in $[m_k]=\{1,\dots,m_k\}$.  Consider a model for the random vector
  $(X_1,\dots,X_p)$ that has the ideal of model invariants $I \subset
  \bbr[p_{i_1, \ldots, i_p}]$.  Let $J \subset \bbr[q_{i_1, \ldots,
    i_{p-1}}, p_{i_1, \ldots, i_p}]$ be the ideal
  $$J \quad = \quad I \, \, \, + \, \, \, \bigg\langle q_{i_1,\ldots, i_{p-1}} -
    \sum_{j =1}^{m_p} p_{i_1, \ldots, i_{p-1} j} \, \, | \, \, i_k \in
    [m_k] \bigg\rangle .$$
  Then the elimination ideal $J \cap \bbr[q_{i_1,
    \ldots, i_{p-1}}]$ is the ideal of model invariants of the model
  created by marginalizing to $X' = (X_1, \ldots, X_{p-1})$.
\end{proposition}

Up to this point, we have made very little use of the inequality
constraints that can arise in the definition of a semi-algebraic set.  In
both of our conditional independence models, the inequality constraints
arose from the fact that we needed to generate a probability distribution,
and were supplied by the positive definite cone or the probability simplex.
In general, however, we may need non-trivial inequality constraints to
describe the model.  Currently, very little is known about the needed
inequality constraints, even in simple examples.  This occurs, for
instance, in the marginalization of conditional independence models.

\begin{example}[Marginalization of an Independence Model]\rm
  Let $A$ be the semi-algebraic set of probability vectors for a discrete
  random vector $X = (X_1, X_2, X_3)$ satisfying the conditional
  independence constraint $X_1 \ind X_2 \mid X_3$.  Let $\psi(A)$ denote
  the image of this model after marginalizing out the random variable
  $X_3$.
  
  The joint distribution of $X_1$ and $X_2$ can be represented as a matrix
  $(p_{ij})$.  Assuming as above that $X_k$ takes on values in $[m_k]$, the
  conditional independence constraint $X_1 \ind X_2 \mid X_3$ implies that the
  matrix $(p_{ij})$ has rank less than or equal to $m_3$.  The set of
  equality constraints that arise from this parametrization are the set of
  $(m_3+1) \times (m_3 + 1)$ minors of the matrix $(p_{ij})$.  However, it
  is not true that these equality constraints together with the inequality
  constraints arising from the probability simplex suffice to define this
  model.  The smallest example of this occurs when, $m_1 = m_2 = 4$ and
  $m_3 = 3$.  In this case the ideal $I(\psi(A))$ is generated by the
  determinant of the generic $4 \times 4$ matrix $(p_{ij})$.  Fix a small
  value of $\epsilon > 0$. The matrix
  $$\frac{1}{8(1 + \epsilon)}\begin{pmatrix}
    1 & 1 & \epsilon & \epsilon \\
    \epsilon & 1 & 1 & \epsilon \\
    \epsilon & \epsilon & 1  & 1 \\
    1 & \epsilon & \epsilon & 1 \end{pmatrix}$$
  represents a probability
  distribution that satisfies the determinant constraint (the matrix has
  rank 3). However, it can be shown that this probability distribution does
  not belong to $\psi(A)$.  That is, this bivariate distribution is not the
  marginalization of a trivariate distribution exhibiting conditional
  independence.  Thus, in addition to the equality constraint, there are
  non-trivial inequality constraints that define the marginalized
  independence model.  More about this example can be found in
  \cite{Mond2003}.  \qed
\end{example}
\mbox{ }


\par

\setcounter{section}{5}
\setcounter{equation}{0} 
\noindent {\bf 5. Solving likelihood equations}

Let $\mathcal{P}=(P_\eta\mid \eta\in N)$ be a regular exponential family
with canonical sufficient statistic $T$.  If we draw a sample $X_1,\dots,X_n$ of
independent random vectors from $P_\eta$, then, as detailed in Section 2,
the canonical statistic becomes $\sum_{i=1}^n T(X_i) =: n\bar T$ and the
log-likelihood function takes the form
\begin{equation}
  \label{eq:loglikfctsample}
  \ell(\eta\mid \bar T)=n\,[ \eta^t \bar T - 
  \phi(\eta)]
\end{equation}
For maximum likelihood estimation in an algebraic exponential family
$\mathcal{P}_M=(P_\eta\mid \eta\in M)$, $M\subseteq N$, we need to maximize
$\ell(\eta\mid\bar T)$ over the set $M$.  

Let $A$ and $g$ be the semi-algebraic set and the diffeomorphism that
define the parameter space $M$.  Let $I(A)=\langle f_1,\dots,f_m\rangle$ be
the ideal of model invariants and $\gamma=g(\eta)$ the parameters after
reparametrization based on $g$.  If boundary issues are of no concern then
the maximization problem can be relaxed to 
\begin{equation}
\label{eq:mlproblem}
\begin{array}{l}
\max \; \ell(\gamma\mid \bar T) \\
\text{subject to} \; f_i(\gamma) = 0,\quad i=1,\dots,m, 
\end{array}
\end{equation}
where
\begin{equation}
  \label{eq:ellgamma}
  \ell(\gamma\mid \bar T) =g^{-1}(\gamma)^t \bar T - \phi(g^{-1}(\gamma)).  
\end{equation}
If $\ell(\gamma\mid \bar T)$ has rational partial derivatives then the
maximization problem (\ref{eq:mlproblem}) can be solved algebraically by
solving a polynomial system of critical equations.  Details on this
approach in the case of discrete data can be found in
\cite{catanese:2006,hosten:2005}.  However, depending on the interplay of
$g^{-1}$ and the mean parametrization $\zeta$, which according to
(\ref{eq:meanmap}) is the gradient map of the log-Laplace transform $\phi$,
such an algebraic approach to maximum likelihood estimation is possible
also in other algebraic exponential families.

\begin{proposition}
  \label{prop:lagrangealgebraic}
  The function $\ell(\gamma\mid \bar T)$ has rational partial derivatives
  if (i) the map $\zeta\circ g^{-1}$ is a rational map and (ii) the map
  $g^{-1}$ has partial derivatives that are rational functions.
\end{proposition}

\begin{example}[Discrete likelihood equations]
  \rm For the discrete exponential family from Example \ref{ex:discrete},
  the mean parameters are the probabilities $p_1,\dots,p_{m-1}$. The
  inverse of the mean parametrization map has component functions
  $(\zeta^{-1})_x=\log(p_x/p_m)$, where $p_m=1-p_1-\dots-p_{m-1}$.  Since
  $d\log(t)/dt=1/t$ is rational, $\zeta^{-1}$ has rational partial
  derivatives.  Hence, maximum likelihood estimates can be computed
  algebraically if the discrete algebraic exponential family is defined in
  terms of the probability coordinates $p_1,\dots,p_{m-1}$.  This is the
  context of the above mentioned work by \cite{catanese:2006,hosten:2005}.
  \qed
\end{example}

\begin{example}[Factor analysis]
  \label{ex:gauss-mle}
  \rm The mean parametrization $\zeta$ for the family of multivariate
  normal distributions and its inverse $\zeta^{-1}$ are based on matrix
  inversions and thus are rational maps.  Thus algebraic maximum likelihood
  estimation is possible whenever a Gaussian algebraic exponential family
  is defined in terms of coordinates $g(\eta)$ for a rational map $g$.
  This includes families defined in the mean parameters $(\mu,\Sigma)$ or
  the natural parameters $(\Sigma^{-1}\mu,\Sigma^{-1})$.
  
  As a concrete example, consider the factor analysis model with one factor
  and four observed variables.  In centered form this model is the family
  of multivariate normal distributions $\ND_4(0,\Sigma)$ on $\RRR^4$ with
  positive definite covariance matrix
  \begin{equation}
    \label{eq:F}
      \Sigma = \text{diag}(\omega) + \lambda\lambda^t,
  \end{equation}
  where $\omega\in (0,\infty)^4$ and $\lambda\in\RRR^4$.  Equation
  (\ref{eq:F}) involves polynomial expressions in
  $\theta=(\omega,\lambda)$.  For algebraic maximum likelihood estimation,
  however, it is computationally more efficient to employ the fact that condition
  (\ref{eq:F}) is equivalent to requiring that the positive definite
  natural parameter $\Sigma^{-1}$ can be expressed as
  \begin{equation}
    \label{eq:Finv}
    \Sigma^{-1}(\theta) = \text{diag}(\omega) -
    \lambda\lambda^t,
  \end{equation}
  with $\theta=(\omega,\lambda)\in(0,\infty)^4\times\RRR^4$; compare
  \citet[\S8]{drton:2006}.  When parametrizing $\Sigma^{-1}$ the map $g$ is
  the identity map.  
  
  Let $S$ be the empirical covariance matrix from a sample of random
  vectors $X_1,\dots,X_n$ in $\RRR^4$; compare (\ref{eq:empiricalXS}).  We
  can solve the maximization problem (\ref{eq:mlproblem}) by plugging the
  polynomial parametric expression for $\gamma=\Sigma^{-1}$ from
  (\ref{eq:Finv}) into the Gaussian version of the log-likelihood function
  in (\ref{eq:ellgamma}).  Taking partial derivatives we find the equations
  \begin{equation}
    \label{eq:fa-deriv}
    \frac{1}{\det(\Sigma^{-1}(\theta))}\cdot\frac{
      \partial \det(\Sigma^{-1}(\theta))}{
      \partial \theta_i} =\trace\bigg[S\cdot \frac{
      \partial \Sigma^{-1}(\theta)}{
      \partial \theta_i}\bigg], \qquad i=1,\dots,8.
  \end{equation}
  These equations can be made polynomial by multiplying by
  $\det(\Sigma^{-1}(\theta))$.  Clearing the denominator introduces many
  additional solutions $\theta\in\mathbb{C}^8$ to the system, which lead to
  non-invertible matrices $\Sigma^{-1}(\theta)$.  However, these extraneous
  solutions can be removed using an operation called {\em saturation\/}.
  After saturation, the (complex) solution set of (\ref{eq:fa-deriv}) is
  seen to consist of $57$ isolated points.  These 57 solutions come in
  pairs $\theta_\pm=(\omega,\pm\lambda)$; one solution has $\lambda=0$.
  
  When the empirical covariance matrix $S$ is rounded then we can compute
  the 57 solutions using software for algebraic and numerical solving of
  polynomial equations.  For the example
  \[
  S_1=
  \begin{pmatrix}
    13&2&-1&3\\
    2&11&3&2\\
    -1&3&9&1\\
    3&2&1&7
  \end{pmatrix}
  \]
  we find that (\ref{eq:fa-deriv}) has 11 feasible solutions in
  $(0,\infty)^4\times\RRR^4$.  Via (\ref{eq:Finv}), these solutions define
  6 distinct factor analysis covariance matrices.  Two of these matrices
  yield local maxima of the likelihood function:
  \[
  \begin{pmatrix}
13 &  2.1242 &  0.9870 & 2.5876\\
2.1242 & 11 & 0.89407 & 2.3440\\
0.9870 & 0.8941 & 9 & 1.0891\\
2.5876 & 2.3440 & 1.0891 &7
  \end{pmatrix},
  \qquad
  \begin{pmatrix}
 13 &  2.1816 & 1.0100 & 1.0962\\
  2.1816 & 11 & 2.3862 & 2.3779\\
  1.0100 & 2.3862 & 9 & 1.1990\\
  1.0962 &  2.3779 & 1.1990 & 7    
  \end{pmatrix}.
  \]
  The matrix to the left has the larger value of the likelihood function
  and we claim that it yields the global maximum.  For this claim to be
  valid we have to check that no matrix close to the boundary of the set
  $\{\Sigma^{-1}(\theta)\mid\theta\in(0,\infty)^4\times \RRR^4\}$ has
  larger value of the likelihood function.  Suppose this was not true.
  Then the likelihood function would have to achieve its global maximum
  over the cone of positive definite matrices outside the set
  $\{\Sigma^{-1}(\theta)\mid\theta\in(0,\infty)^4\times \RRR^4\}$.  In
  order to rule out this possibility, we consider all the complex solutions
  $\theta\not\in (0,\infty)^4\times \RRR^4$ of (\ref{eq:fa-deriv}) that
  induce real and positive definite matrices $\Sigma^{-1}(\theta)$.  There
  are ten such solutions, which all have $\omega\in\RRR^4$ and purely
  imaginary $\lambda\in i\RRR^4$.  There are five different induced
  matrices $\Sigma^{-1}(\theta)$, but at all of them the likelihood
  function is smaller than for the two quoted local maximizer.  This
  confirms our claim.
  
  As a second interesting example consider
  \[
  S_2=
  \begin{pmatrix}
    31& 11& -1& 5\\
    11& 23& 3& -2\\
    -1& 3& 7& 1\\
    5& -2& 1& 7
  \end{pmatrix}.
  \]
  The equations (\ref{eq:fa-deriv}) have again 11 feasible solutions
  $\hat\theta$.  Associated are 6 distinct factor analysis covariance
  matrices that all correspond to saddle points of the
  likelihood function.  Hence, if we close the set of inverse covariance
  matrices $\{\Sigma^{-1}(\theta)\mid\theta\in(0,\infty)^4\times \RRR^4\}$,
  then the global optimum of the likelihood function over this closure must
  be attained on the boundary.
  
  In order to determine which boundary solution provides the global maximum
  of the likelihood function, it is more convenient to switch back to the
  standard parameterization in (\ref{eq:F}), which writes the covariance
  matrix as $\Sigma(\theta)$ for $\theta=(\omega,\lambda)$ in
  $(0,\infty)^4\times \RRR^4$.  The closure of $\{\Sigma(\theta)\mid
  \theta\in(0,\infty)^4\times\RRR^4\}$ is obtained by closing the parameter
  domain to $[0,\infty)^4\times\RRR^4$.  Since $S_2$ is positive definite,
  the global maximizer of the likelihood function must be a matrix of full
  rank, which implies that at most one of the four parameters $\omega_i$
  can be zero.  In each of the four possible classes of boundary cases the
  induced likelihood equations (in 7 parameters) have a closed form
  solution leading to a unique covariance matrix.  We find that the global
  maximum is achieved in the case $\omega_1=0$.  The global maximizer of
  the likelihood functions over the closure of the parameter space equals
  \[
  \begin{pmatrix}
     31 & 11 & -1 & 5\\
     11 & 23 & -0.3548 & 1.7742\\
     -1&-0.3548&7&-0.1613\\
     5&1.7742&-0.1613&7
  \end{pmatrix}.
  \]
  In the factor analysis literature data leading to such boundary problems
  are known as Heywood cases.  Hence, our computation {\em proves\/} that
  $S_2$ constitutes a Heywood case.  \qed
\end{example}
\mbox{ }
\par

\noindent {\bf 6.  Conclusion}

In this paper, we have attempted to present a useful, unified definition of
an algebraic statistical model.  In this definition, an algebraic model is
a submodel of a reference model with nice statistical properties.  Working
primarily with small examples of conditional independence models, we have
tried to illustrate how our definition might be a useful framework, in
which the geometry of parameter spaces can be related to properties of
statistical inference procedures.  Since we impose algebraic structure,
this geometry can be studied using algebraic techniques, which
allow one to tackle problems where simple linear arguments will not work.
In order to apply these algebraic techniques in a particular example of
interest, one can resort to one of the many software systems, both free and
commercial, that provide implementations of algorithms for carrying out the
necessary computations.  A comprehensive list of useful software can be
found in Chapter 2 of \cite{ascb}.

While we believe that future work in algebraic statistics may involve
reference models in which the notion of ``nice statistical properties'' is
filled with life in many different ways, we also believe that the most
important class of reference models are regular exponential families.  This
led us to consider what we termed algebraic exponential families.  These
families were shown to be flexible enough to encompass structures arising
from marginalization, i.e., the involvement of hidden variables.  Hidden
variable models typically do not form curved exponential families, which
triggered \cite{geiger:2001} to introduce their stratified exponential
families. These stratified families are more general than both algebraic
and curved exponential families but, as our Example \ref{ex:stratified}
suggests, they seem in fact to be too general to allow the derivation of
results that would hold in the entire class of models.  In algebraic
exponential families, on the other hand, the restriction to semi-algebraic
sets entails that parameter spaces always have nice local geometric
properties and phenomena as created in Example \ref{ex:stratified} cannot
occur.  In light of this fact, our algebraic exponential families appear to
be in particular a good framework for the study of hidden variable models,
which are widely used models whose statistical properties have yet to be
understood in entirety.



\bigskip

\noindent {\large\bf Acknowledgment}

Mathias Drton was partially supported by the NSF (DMS-0505612).  


\def\refname{\large\bf References}


\vskip .65cm
\noindent
Department of Statistics, University of Chicago, Chicago, IL, U.S.A.
\vskip 2pt
\noindent
E-mail: (drton@galton.uchicago.edu)
\vskip 2pt
\noindent
Society of Fellows, Harvard
  University, Cambridge, MA, U.S.A.
\vskip 2pt
\noindent
E-mail: (seths@math.harvard.edu)
\vskip .3cm
\end{document}